\documentclass[12pt]{article}
\usepackage{amssymb,amsmath}
\usepackage{graphicx}

\newenvironment{keywords}{\centerline{\bf\small
Keywords}\begin{quote}\small}{\par\end{quote}\vskip 1ex}

\newtheorem{thm}{Theorem}

\newtheorem{exmp}{Example}
\newtheorem{cor}{Corollary}
\newtheorem{rem}{Remark}

\newcommand{\eoe}{\hspace*{\fill} $\diamondsuit\quad$}

\renewcommand{\P}{{\rm P}}
\newcommand{\E}{{\bf E}}
\newcommand{\uP}{{\overline\P}}
\newcommand{\uE}{{\overline\E}}
\newcommand{\lP}{{\underline\P}}
\newcommand{\lE}{{\underline\E}}

\newcommand{\erre}{\mathbf{R}}
\newcommand{\enne}{\mathbf{N}}
\newcommand{\giv}{{\,|\,}}

\newcommand{\ucat}{{\mathcal{X}=\{x_1,\ldots ,x_k\}}}

\newcommand{\uc}{{\mathcal{X}}}

\newcommand{\ill}{I}
\newcommand{\healthy}{H}
\newcommand{\postest}{+}
\newcommand{\negtest}{-}

\newcommand{\ux}{{\mathrm{X}}}

\newcommand{\us}{{\mathrm{S}}}

\newcommand{\uxs}{{\mathbf{x}}}
\newcommand{\uys}{{\mathbf{y}}}
\newcommand{\uds}{{\mathbf{d^i}}}
\newcommand{\uesi}{{\mathbf{e^i}}}
\newcommand{\uss}{{\mathbf{s}}}

\newcommand{\ufs}{{\mathbf{f}}}
\newcommand{\latent}{{\mathbf{X}}}
\newcommand{\latentfut}{{\mathbf{X'}}}
\newcommand{\manifest}{{\mathbf{S}}}
\newcommand{\realized}{{\mathbf{X}}}

\newcommand{\ut}{{\mathbf{t}}}
\newcommand{\teta}{{\mathbf{\theta}}}
\newcommand{\ste}{{\Theta}}

\newcommand{\ssette}{{\Theta:=\{\mathbf{\theta}=(\theta_1,\ldots,\theta_k)\,|\, \sum_{i=1}^k \theta
_i=1,\; 0\leq \theta_i\leq 1\}}}
\newcommand{\sett}{{\mathcal{T}}}
\newcommand{\ssett}{{\mathcal{T}:=\{\ut=(t_1,\ldots ,t_k)\,|\,
\sum_{j=1}^k\, t_k=1,\, 0<t_j<1\}}}

\newcommand{\seqpn}{{(p_n)_{n\in\enne}}}

\newcommand{\mzero}{{\mathcal{M}_0}}

\begin{document}

\title{\vspace{-4ex}
\normalsize\sc Technical Report \hfill IDSIA-05-07
\vskip 2mm\bf\Large\hrule height5pt \vskip 4mm
Learning about a Categorical Latent Variable under Prior Near-Ignorance
\vskip 4mm \hrule height2pt}
\author{
\normalsize Alberto Piatti, IDSIA, Switzerland, alberto.piatti@idsia.ch \\
\normalsize Marco Zaffalon, IDSIA, Switzerland, zaffalon@idsia.ch \\
\normalsize Fabio Trojani, University of S. Gallen, Switzerland, fabio.trojani@unisg.ch \\
\normalsize Marcus Hutter, ANU\&NICTA, Australia ,marcus@hutter1.net
}
\date{May 2007}

\maketitle

\begin{abstract}
It is well known that complete prior ignorance is not compatible
with learning, at least in a coherent theory of (epistemic)
uncertainty. What is less widely known, is that there is a state
similar to full ignorance, that Walley calls \emph{near-ignorance},
that permits learning to take place. In this paper we provide new
and substantial evidence that also near-ignorance cannot be really
regarded as a way out of the problem of starting statistical
inference in conditions of very weak beliefs. The key to this result
is focusing on a setting characterized by a variable of interest
that is \emph{latent}. We argue that such a setting is by far the
most common case in practice, and we show, for the case of
categorical latent variables (and general \emph{manifest} variables)
that there is a sufficient condition that, if satisfied, prevents
learning to take place under prior near-ignorance. This condition is
shown to be easily satisfied in the most common statistical
problems.
\def\contentsname{\centering\normalsize Contents}
{\parskip=-2.7ex\tableofcontents}
\end{abstract}


\section{Introduction}\label{sec:introduction}

Epistemic theories of statistics are often concerned with the
question of \emph{prior ignorance}. Prior ignorance means that a
subject, who is about to perform a statistical analysis, has not any
substantial belief about the underlying data-generating process.
Yet, the subject would like to exploit the available sample to draw
some statistical inference, i.e., the subject would like to use the
data to learn, moving away from the initial condition of ignorance.
This situation is very important as it is often desirable to start a
statistical analysis with weak assumptions about the problem of
interest, thus trying to implement an objective-minded approach to
statistics.

A fundamental question is if prior ignorance is compatible with
learning. Walley gives a negative answer for the case of his
self-consistent (or \emph{coherent}) theory of statistics: he shows,
in a very general sense, that \emph{vacuous} prior beliefs lead to
vacuous posterior beliefs, irrespective of the type and amount of
observed data \cite[Section~7.3.7]{Walley1991}. But, at the same
time, he proposes focusing on a slighlty different state of beliefs,
called \emph{near-ignorance}, that does enable learning to take
place \cite[Section~4.6.9]{Walley1991}. Loosely speaking,
near-ignorant beliefs are beliefs close but not equal to vacuous
(see Section~\ref{sec:nearignorance}). The possibility to learn
under prior near-ignorance is shown, for instance, in the special
case of the near-ignorance prior defining the \emph{imprecise
Dirichlet model} (IDM). This is a popular model used in the case of
inference from categorical data generated by a discrete process
(\cite{Walley1996, Bernard2005}).

In this paper, we also focus on a categorical random variable $\ux$,
expressing the outcomes of a multinomial process, but we assume that
such a variable is \emph{latent}. This means that we cannot observe
the realizations of $\ux$, so we can learn about it only by means of
another (not necessarily categorical) variable $\us$, related to
$\ux$ in some known way. Variable $\us$ is assumed to be
\emph{manifest}, in the sense that its realizations can be observed
(see Section~\ref{sec:latent}).

In such a setting, we introduce a condition in
Section~\ref{sec:conditions}, related to the likelihood of the
observed data, that is shown to be sufficient to prevent learning
about $\ux$ under prior near-ignorance. The condition is very
general as it is developed for any prior that models near-ignorance
(not only the one used in the IDM), and for very general kinds of
relation between $\ux$ and $\us$. We show then, by simple examples,
that such a condition is easily satisfied, even in the most
elementary and common statistical problems.

In order to appreciate this result, it is important to realize that
latent variables are ubiquitous in problems of uncertainty. It can
be argued, indeed, that there is a persistent distinction between
(latent) facts (e.g., health, state of economy, color of a ball) and
(manifest) observations of facts: one can regard them as being
related by a so-called \emph{observational process}; and the point
is that these kinds of processes are imperfect in practice.
Observational processes are often neglected in statistics, when
their imperfection is deemed to be tiny. But a striking outcome of
the present research is that, no matter how tiny the imperfection,
provided it exists, learning is not possible under prior
near-ignorance.

In our view, the present results raise serious doubts about the
possibility to adopt a condition of prior near-ignorance in real, as
opposed to idealized, applications of statistics.  As a consequence,
it may make sense to consider re-focusing the research about this
subject on developing models of very weak states of belief that are,
however, stronger than near-ignorance.

\section{Categorical Latent Variables}\label{sec:latent}

In this paper, we follow the general definition of
\emph{latent} and \emph{manifest variables} given by
\cite{Skrondal2004}: a \emph{latent variable} is a random variable
whose realizations are unobservable (hidden), while a \emph{manifest
variable} is a random variable whose realizations can be directly
observed. The concept of latent variable is central in many
sciences, like for example psychology and medicine.
\cite{Skrondal2004} list several fields of application and several
phenomena that can be modeled using latent variables, and conclude
that latent variable modeling ``\emph{pervades modern mainstream
statistics},'' although ``\emph{this omni-presence of latent
variables is commonly not recognized, perhaps because latent
variables are given different names in different literatures, such
as random effects, common factors and latent classes},'' or hidden variables.

But what are latent variables in practice? According to
\cite{Boorsbom2002}, there may be different interpretations of
latent variables. A latent variable can be regarded, for example, as
an unobservable random variable that exists independently of the
observation. An example is the unobservable health status of a
patient that is subject to a medical test. Another possibility is to
regard a latent variable as a product of the human mind, a construct
that does not exist independent of the observation. For example the
\emph{unobservable state of the economy}, often used in economic
models. In this paper, we assume the existence of a latent
categorical random variable $\ux$, with outcomes in $\ucat$ and
unknown chances $\teta\in\ssette$, without stressing any particular
interpretation.

Suppose now that our aim is to predict, after $N$ realizations of
the variable $\ux$, the next outcome (or the next $N'$ outcomes).
Because the variable $\ux$ is latent and therefore unobservable by
definition, the only possible way to learn something about the
probabilities of the next outcome is to observe the realizations of
some manifest variable $\us$ related, in a known way, to the
(unobservable) realizations of $\ux$. An example of known
relationship between latent and manifest variables is the following.

\begin{exmp}\rm
\label{exmp:introtest} We consider a binary medical diagnostic test
used to assess the health status of a patient with respect to a
given disease. The accuracy of a diagnostic test\footnote{For
further details about the modeling of diagnostic accuracy with
latent variables see \cite{Yang1997}.} is determined by two
probabilities: the \emph{sensitivity} of a test is the probability
of obtaining a positive result if the patient is diseased; the
\emph{specificity} is the probability of obtaining a negative result
if the patient is healthy. Medical tests are assumed to be imperfect
indicators of the unobservable true disease status of the patient.
Therefore, we assume that the probability of obtaining a positive
result when the patient is healthy, respectively of obtaining a
negative result if the patient is diseased, are non-zero. Suppose,
to make things simpler, that the sensitivity and the specificity of
the test are known. In this example, the unobservable health status
of the patient can be considered as a binary latent variable $\ux$
with values in the set $\{\textrm{Healthy},\textrm{Ill}\}$, while
the result of the test can be considered as a binary manifest
variable $\us$ with values in the set $\{\textrm{Negative
result},\textrm{Positive result}\}$. Because the sensitivity and the
specificity of the test are known, we know how $\ux$ and $\us$ are
related.\eoe
\end{exmp}

We continue discussion about this example later on, in the light of
our results, in Example \ref{exmp:diagnostictest} of Section
\ref{sec:conditions}.

\section{Near-Ignorance Priors}\label{sec:nearignorance}

Consider a categorical random variable $\ux$ with outcomes in
$\ucat$ and unknown chances $\teta\in\ste$. Suppose that we have no
relevant prior information about $\teta$ and we are therefore in a
situation of prior ignorance. How should we model our prior beliefs
in order to reflect the initial lack of knowledge?

Let us give a brief overview of this topic in the case of coherent
models of uncertainty, such as Bayesian probability and Walley's
theory of \emph{coherent lower previsions}.

In the traditional Bayesian setting, prior beliefs are modeled
using a single prior probability distribution. The problem of
defining a standard prior probability distribution modeling a
situation of prior ignorance, a so-called \emph{noninformative
prior}, has been an important research topic in the last two
centuries\footnote{Starting from the work of Laplace at the
beginning of the 19$^\text{th}$ century (\cite{Laplace1820}).} and, despite
the numerous contributions, it remains an open research issue, as
illustrated by \cite{Kass1996}. See also \cite{Hutter2006} for recent
developments and complementary considerations. There are many
principles and properties that are desirable to model a situation of
prior ignorance and that have been used in past research to define
noninformative priors. For example Laplace's \emph{symmetry or
indifference} principle has suggested, in case of finite possibility
spaces, the use of the uniform distribution. Other principles, like
for example the principle of \emph{invariance under group transformations},
the \emph{maximum entropy} principle, the \emph{conjugate
priors} principle, etc., have suggested the use of other
noninformative priors, in particular for continuous possibility
spaces, satisfying one or more of these principles. But, in general,
it has proven to be difficult to define a standard noninformative
prior satisfying, at the same time, all the desirable principles.

In the case of finite possibility spaces, we agree with
\cite{DeCooman2006} when they say that there are at least two
principles that should be satisfied to model a situation of prior
ignorance: the \emph{symmetry principle} and the \emph{embedding
principle}. The \emph{symmetry principle} states that, if we are
completely ignorant a priori about $\teta$, then we have no reason
to favour one possible outcome of $\ux$ to another, and therefore
our probability model on $\theta$ should be symmetric. This
principle recalls Laplace's \emph{symmetry or indifference}
principle that, in the past decades, has suggested the use of the
\emph{uniform prior} as standard noninformative prior. The
\emph{embedding principle} states that, for each possible event $A$,
the probability assigned to $A$ should not depend on the possibility
space $\uc$ in which $A$ is embedded. In particular, the probability
assigned a priori to the event $A$ should be invariant with respect
to refinements and coarsenings of $\uc$. It is easy to show that the
embedding principle is not satisfied by the uniform distribution.
How should we model our prior ignorance in order to satisfy these
two principles? \cite{Walley1991} gives a compelling answer to this
question: he proves\footnote{In Note~7, p.~526. See also
Section~5.5.} that the only probability model consistent with
coherence and with the two principles is the \emph{vacuous
probability model}, i.e., the model that assigns, for each
non-trivial event $A$, lower probability $\lP(A)=0$ and upper
probability $\uP(A)=1$. It is evident that this model cannot be
expressed using a single probability distribution. It follows that,
to model properly and in a coherent way a situation of prior
ignorance, we need \emph{imprecise probabilities}.\footnote{For a
complementary point of view, see \cite{Hutter2006}.}

Unfortunately, adopting the vacuous probability model for $\ux$ is not a
practical solution to our initial problem, because it produces only
vacuous posterior probabilities. \cite{Walley1991} suggests, as
practical solution, the use of \emph{near-ignorance priors}. A near-ignorance prior is a large closed convex set $\mzero$ of probability
distributions for $\theta$, very close to the vacuous probability model, which
produces a priori \emph{vacuous expectation}s for various functions
$f$ on $\ste$, i.e., such that $\lE(f)=\inf_{\theta\in\Theta}
f(\theta)$ and $\uE(f)=\sup_{\theta\in\Theta} f(\theta)$.

An example of near-ignorance prior that is particularly instructive
is the set of priors $\mzero$ used in the \emph{imprecise Dirichlet
model} (IDM). The IDM models a situation of prior ignorance about
the chances $\theta$ of a categorical random variable $\ux$. The
near-ignorance prior $\mzero$ used in the IDM consists in the set of
all Dirichlet densities $p(\theta)=dir_{s,\ut}(\theta)$ for a fixed
$s>0$ and all $\ut\in\sett$, where
\begin{equation}\label{eq:dir}
  dir_{s,\ut}(\theta):=\frac{\Gamma (s)}{\prod_{i=1}^k \Gamma (st_i)}\,\prod_{i=1}^k\, \theta_i^{st_i-1},
\end{equation}
and
\begin{equation}\label{eq:setT}
  \ssett.
\end{equation}
The particular choice of $\mzero$ in the IDM implies vacuous prior
expectations for all functions $f(\theta)=\theta_i^{N'}$, for all
$N'\geq 1$ and all $i\in\{1,\ldots ,k\}$, i.e.,
$\lE(\theta_i^{N'})=0$ and
$\uE(\theta_i^{N'})=1$. Choosing $N'=1$, we have, a priori,
$$
  \lP(\ux=x_i)=\lE(\theta_i)=0,\quad
  \uP(\ux=x_i)=\uE(\theta_i)=1.
$$
It follows that the particular near-ignorance prior $\mzero$ used in
the IDM implies vacuous prior probabilities for each possible
outcome of the variable $\ux$. It can be shown that this particular
set of priors satisfies both the symmetry and embedding principles.

But what is the difference between the vacuous probability model and
the the near-ignorance prior used in the IDM? In fact, although both
models produce vacuous prior probabilities and both models satisfy
the symmetry and embedding principles, the IDM yields posterior probabilities that are not vacuous, while the vacuous probability model produces only
vacuous posterior probabilities. The answer to this question is the
reason why we use the term \emph{near-ignorance}: in the IDM,
although we are completely ignorant about the possible outcomes of
the variable $\ux$, we are not completely ignorant about the chances
$\theta$, because we assume a particular class of prior
distributions, i.e., the Dirichlet distributions for a fixed value
of $s$.

\section{Limits of Learning under Prior Near-Ignorance}\label{sec:conditions}

Consider a sequence of independent and identically distributed (IID)
categorical latent variables $(\ux_i)_{i\in\enne}$ with outcomes in
$\uc$ and unknown chances $\theta\in\ste$, and a sequence of
independent manifest variables $(\us_i)_{i\in\enne}$. We assume that
a realization of the manifest variable $\us_i$ can be observed only
after an (unobservable) realization of the latent variable $\ux_i$
and that the probability distribution of $\us_i$ given $\ux_i$ is
known for each $i\in\enne$. Furthermore, we assume $\us_i$ to be
independent of the chances $\theta$ of $\ux_i$ given $\ux_i$. Define
the random variables $\latent:=(\ux_1,\ldots,\ux_N)$,
$\manifest:=(\us_1,\ldots,\us_N)$ and
$\latentfut:=(\ux_{N+1},\ldots,\ux_{N+N'})$.

We focus on the problem of predictive inference.\footnote{For a
general presentation of predictive inference see \cite{Geisser1993};
for a discussion of the imprecise probability approach to predictive
inference see \cite{Walley1999}.} Suppose that we observe a dataset
$\uss$ of realizations of manifest variables $\us_1,\ldots,\us_N$
related to the (unobservable) dataset $\uxs\in\uc^N$ of realizations
of the variables $\ux_1,\ldots,\ux_N$. Using the notation defined
above we have $\manifest=\uss$ and $\latent=\uxs$. Our aim is to
predict the outcomes of the next $N'$ variables
$\ux_{N+1},\ldots,\ux_{N+N'}$. In particular, given
$\uxs'\in\uc^{N'}$, our aim is to calculate
$\lP(\latentfut=\uxs'\giv\manifest=\uss)$ and
$\uP(\latentfut=\uxs'\giv\manifest=\uss)$. To simplify notation,
when no confusion is possible, we denote in the rest of the paper
$\manifest=\uss$ with $\uss$ and $\latentfut=\uxs'$ with $\uxs'$.
The (in)dependence structure can be depicted graphically as follows:

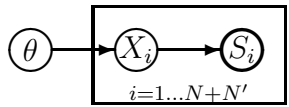
\begin{center}
\unitlength=0.8pt
\begin{picture}(130,45)(10,5)
\thicklines \put(20,30){\circle{20}\makebox(0,0)[cc]{$\theta$}}
\put(70,30){\circle{20}\makebox(0,0)[cc]{$X_i$}}
\put(120,30){\circle{19}}
\put(120,30){\circle{21}\makebox(0,0)[cc]{$S_i$}}
\put(50,5){\framebox(90,45)[cb]{$\scriptstyle i=1\ldots N+N'$}}
\put(30,30){\vector(1,0){30}} \put(80,30){\vector(1,0){30}}
\end{picture}
\end{center}

Modelling our prior ignorance about the parameters $\theta$ with a
near-ignorance prior $\mzero$ and denoting by
$\mathbf{n'}:=(n_1',\ldots,n_k')$ the frequencies of the dataset
$\uxs'$, we have
\begin{eqnarray}\nonumber
\lP(\mathbf{\uxs'}\giv\uss) & = &
\inf_{p\in\mzero}\P_p(\mathbf{\uxs'}\giv\uss):=\\\nonumber & = &
\inf_{p\in\mzero}\int_{\ste} \prod_{i=1}^k \theta_i^{n_i'}
p(\theta\giv\uss) d\theta=\\\nonumber & =: &
\inf_{p\in\mzero}\E_p\left(\prod_{i=1}^k
\theta_i^{n_i'}\giv\uss\right) =\\\nonumber & = &
\lE\left(\prod_{i=1}^k
\theta_i^{n_i'}\giv\uss\right),\label{eq:lowerexpected}
\end{eqnarray}
where, according to Bayes theorem, $$p(\theta\giv\uss)=
\frac{\P(\uss\giv\theta) p(\theta)}{\int_{\ste} \P(\uss\giv\theta)
p(\theta) d\theta},$$ provided that $\int_{\ste} \P(\uss\giv\theta)
p(\theta) d\theta\neq 0$. Analogously, substituting $\sup$ to $\inf$
in (\ref{eq:lowerexpected}), we obtain
\begin{equation}\label{eq:upperexpected}
\uP(\mathbf{\uxs'}\giv\uss)=\uE\left(\prod_{i=1}^k
\theta_i^{n_i'}\giv\uss\right).
\end{equation}
The central problem now is to choose $\mzero$ so as to be as
ignorant as possible a priori and, at the same time, to be able to
learn something from the observed dataset of manifest variables
$\uss$. Theorem \ref{thm:fondamentale} and the following corollaries
yield a first partial solution to the above problem, stating several
conditions for learning under prior near-ignorance.

\begin{thm}\label{thm:fondamentale}
Let $\uss$ be given. Consider a bounded continuous function $f$
defined on $\Theta$ and denote with $f_{\max}$ the Supremum of $f$
on $\Theta$. If the likelihood function $\P(\uss\giv\theta)$ is
strictly positive\footnote{The Assumption about $\P(\uss\giv\theta)$
in Theorem \ref{thm:fondamentale} can be substituted by the
following weaker assumption. For a given arbitrary small $\delta>0$,
denote with $\Theta_{\delta}$ the measurable set,
$\Theta_{\delta}:=\{\theta\in\Theta\,|\, f(\theta)\geq
f_{\max}-\delta\}.$ If $\P(\uss\giv\theta)$ is such that, $
\lim_{\delta\rightarrow
0}\inf_{\theta\in\Theta_{\delta}}\P(\uss\giv\theta)=c>0, $ then
Theorem \ref{thm:fondamentale} holds.} in each point in which $f$
reaches its maximum value $f_{\max}$ and it is continuous in an
arbitrary small neighborhood of these points, and $\mzero$ is such
that a priori $\uE(f)=f_{\max}$, then
$$\uE(f\giv\uss)=\uE(f)=f_{\max}.$$
\end{thm}

Many corollaries to Theorem \ref{thm:fondamentale} are listed in
Section \ref{app:corollaries} of the Appendix. Here we discuss only
the most important corollary. Consider, given a dataset $\uxs'$, the
particular function $f(\theta)=\prod_{i=1}^k \theta_i^{n_i'}$. This
function is particularly important for predictive inference, because
its lower and upper expectations correspond to the lower and upper
probabilities assigned to the dataset $\uxs'$. It is easy to show
that, in this case, the minimum of $f$ is 0 and is reached in all
the points $\theta\in\ste$ with $\theta_i=0$ for some $i$ such that
$n_i'>0$, while the maximum of $f$ is reached in a single point of
$\ste$ corresponding to the relative frequencies $\mathbf{f'}$ of
the sample $\uxs'$, i.e., at
$\mathbf{f'}=\left(\frac{n_1'}{N'},\ldots,\frac{n_k'}{N'}\right)\in\Theta$,
and the maximum of $f$ is given by $\prod_{i=1}^k
\left(\frac{n_i'}{N'}\right)^{n_i'}$. It follows that vacuous
probabilities regarding the dataset $\uxs'$ are given by
$$\lP(\uxs')=\lE\left(\prod_{i=1}^k \theta_i^{n_i'}\right)=0,$$
$$\uP(\uxs')=\uE\left(\prod_{i=1}^k \theta_i^{n_i'}\right)=\prod_{i=1}^k
\left(\frac{n_i'}{N'}\right)^{n_i'}.$$

\begin{cor}\label{cor:fondamentale4}
Let $\uss$ be given and let $\P(\uss\giv\theta)$ be a continuous
strictly positive function on $\Theta$. Then, if $\mzero$ implies
vacuous prior probabilities for a dataset $\uxs'\in\uc^{N'}$, the
predictive probabilities of $\uxs'$ are vacuous also a posteriori,
after having observed $\uss$, i.e.,
$$\lP(\uxs'\giv\uss)=\lP(\uxs')=0,$$
$$\uP(\uxs'\giv\uss)=\uP(\uxs') =\prod_{i=1}^k
\left(\frac{n_i'}{N'}\right)^{n_i'}.$$
\end{cor}

In other words, Corollary \ref{cor:fondamentale4} states a sufficient
condition that prevents learning to take place under prior
near-ignorance: if the likelihood function $\P(\uss\giv\theta)$ is
continuous and strictly positive on $\Theta$, then all the dataset
$\uxs'\in\uc^{N'}$ for which $\mzero$ implies vacuous probabilities
have vacuous probabilities also a posteriori, after having observed
$\uss$. It follows that, if this sufficient condition is satisfied,
we cannot use near-ignorance priors to model a state of prior
ignorance for the same reason for which, in Section~\ref{sec:nearignorance}, we have excluded the vacuous probability
model: because only vacuous posterior probabilities are produced.

The sufficient condition described above is satisfied very often in
practice, as illustrated by the following striking examples.

\begin{exmp}\label{exmp:diagnostictest}\rm
Consider the medical test introduced in Example \ref{exmp:introtest}
and an (ideally) infinite population of individuals. Denote with the
binary variable $\ux_i\in\{\healthy,\ill\}$ the health status of the
$i$-th individual of the population and with
$\us_i\in\{\postest,\negtest\}$ the results of the diagnostic test
applied to the same individual. We assume that the variables in the
sequence $(\ux_i)_{i\in\enne}$ are IID with unknown chances
$(\theta,1-\theta)$, where $\theta$ corresponds to the (unknown)
proportion of diseased individuals in the population. Denote with
$1-\varepsilon_1$ the sensitivity and with $1-\varepsilon_2$ the
specificity of the test. Then it holds that
$$
  \P(\us_i=\postest\giv\ux_i=\healthy)=\varepsilon_1>0,
$$
$$
  \P(\us_i=\negtest\giv\ux_i=\ill)=\varepsilon_2>0,
$$
where ($\ill,\healthy,\postest,\negtest$) denote
(patient ill, patient healthy, test positive, test negative).

Suppose that we observe the results of the test applied to $N$
different individuals of the population; using our previous notation
we have $\manifest=\uss$. For each individual we have,
\begin{align*}
& \P(\us_i=\postest\giv\theta)=\\
= &
\P(\us_i=\postest\giv\ux_i=\ill)\P(\ux_i=\ill\giv\theta)+\\
+ & \P(\us_i=\postest\giv\ux_i=\healthy)\P(\ux_i=\healthy\giv\theta)=\\
= & \underbrace{(1-\varepsilon_2)}_{>0}\cdot\theta+\underbrace{\varepsilon_1}_{>0}\cdot (1-\theta)>0.\\
\end{align*}
Analogously,
\begin{align*}
& \P(\us_i=\negtest\giv\theta)=\\
= &
\P(\us_i=\negtest\giv\ux_i=\ill)\P(\ux_i=\ill\giv\theta)+\\
+ & \P(\us_i=\negtest\giv\ux_i=\healthy)\P(\ux_i=\healthy\giv\theta)=\\
= & \underbrace{\varepsilon_2}_{>0}\cdot\theta+\underbrace{(1-\varepsilon_1)}_{>0}\cdot (1-\theta)>0.\\
\end{align*}
Denote with $n^{\uss}$ the number of positive tests in the observed
sample $\uss$. Then, because the variables $\us_i$ are independent,
we have
\begin{align*}
&
\P(\manifest=\uss\giv\theta)=((1-\varepsilon_2)\cdot\theta+\varepsilon_1\cdot
(1-\theta))^{n^{\uss}}\cdot\\
\cdot & (\varepsilon_2\cdot\theta+(1-\varepsilon_1)\cdot (1-\theta))^{N-n^{\uss}}>0\\
\end{align*}
for each $\theta\in [0,1]$ and each $\uss\in\uc^N$. Therefore,
according to Corollary \ref{cor:fondamentale4}, all the predictive
probabilities that, according to $\mzero$, are vacuous a priori
remain vacuous a posteriori. It follows that, if we want to avoid
vacuous posterior predictive probabilities, then we cannot model our
prior knowledge (ignorance) using a near-ignorance prior implying
some vacuous prior predictive probabilities. This simple example
shows that our previous theoretical results raise serious questions
about the use of near-ignorance priors also in very simple, common,
and important situations.

The situation presented in this example can be extended, in a
straightforward way, to the general categorical case and has been
studied, in the special case of the near-ignorance prior used in the
imprecise Dirichlet model, in \cite{Piatti2005}.\eoe
\end{exmp}

Example \ref{exmp:diagnostictest} focuses on discrete latent and
manifest variables. In the next example, we show that our
theoretical results have important implications also in models with
discrete latent variables and continuous manifest variables.

\begin{exmp}\label{exmp:continuousmanifest}\rm
Consider the sequence of IID categorical variables
$(\ux_i)_{i\in\enne}$ with outcomes in $\uc^N$ and unknown chances
$\theta\in\Theta$. Suppose that, for each $i\geq 1$, after a
realization of the latent variable $\ux_i$, we can observe a
realization of a continuous manifest variable $\us_i$. Assume that
$p(\us_i\giv\ux_i=x_j)$ is a continuous positive probability
density, e.g., a normal $N(\mu_j,\sigma_j^2)$ density, for each
$x_j\in\uc$. We have
\begin{align*}
p(\us_i\giv\theta) & = \sum_{x_j\in\uc^N} p(\us_i\giv\ux_i=x_j)\cdot
\P(\ux_i=x_j\giv\theta)=\\
& =\sum_{x_j\in\uc^N} \underbrace{p(\us_i\giv\ux_i=x_j)}_{>0}\cdot
\theta_j>0,
\end{align*}
because $\theta_j$ is positive for at least one $j\in\{1,\ldots
,N\}$ and we have assumed $\us_i$ to be independent of $\theta$
given $\ux_i$. Because we have assumed $(\us_i)_{i\in\enne}$ to be a
sequence of independent variables, we have,
$$p(\manifest=\uss\giv\theta)=\prod_{i=1}^N \underbrace{p(\us_i=\uss_i\giv\theta)}_{>0}>0.$$
Therefore, according to Corollary \ref{cor:fondamentale4}, if we
model our prior knowledge using a near-ignorance prior $\mzero$, the
vacuous prior predictive probabilities implied by $\mzero$ remain
vacuous a posteriori. It follows that, if we want to avoid vacuous
posterior predictive probabilities, we cannot model our prior
knowledge using a near-ignorance prior implying some vacuous prior
predictive probabilities.\eoe
\end{exmp}

Examples \ref{exmp:diagnostictest} and \ref{exmp:continuousmanifest}
raise, in general, serious criticisms about the use of
near-ignorance priors in practical applications.

The only predictive model in the literature, of which we are aware,
where a near-ignorance prior is used successfully to obtain
non-vacuous posterior predictive probabilities is the IDM. In the
next example, we explain how the IDM avoids our theoretical
limitations.

\begin{exmp}\label{exmp:idm}\rm
In the IDM, we assume that the IID categorical variables
$(\ux_i)_{i\in\enne}$ are observable. In other words, we have
$\us_i=\ux_i$ for each $i\geq 1$ and therefore the IDM is not a
latent variable model. Having observed $\manifest=\realized=\uxs$,
we have
$$\P(\manifest=\uxs\giv\theta)=\P(\realized=\uxs\giv\theta)=\prod_{i=1}^k \theta_i^{n_i},$$
where $n_i$ denotes the number of times that $x_i\in\uc$ has been
observed in $\uxs$. We have $\P(\realized=\uxs\giv\theta)=0$ for all
$\theta$ such that $\theta_j=0$ for at least one $j$ such that
$n_j>0$ and $\P(\realized=\uxs\giv\theta)>0$ for all the other
$\theta\in\Theta$, in particular for all $\theta$ in the interior of
$\Theta$.

The near-ignorance prior $\mzero$ used in the IDM consists in the
set of all the Dirichlet densities $dir_{s,\ut}(\theta)$ for a fixed
$s>0$ and all $\ut\in\sett$, where
$dir_{s,\ut}(\theta)$ and $\sett$ have been defined in (\ref{eq:dir})
and (\ref{eq:setT}).

The particular choice of $\mzero$ in the IDM implies, for each
$N'\geq 1$ and each $i\in\{1,\ldots ,k\}$, that
$$\lE(\theta_i^{N'})=0,\quad \uE(\theta_i^{N'})=1.$$
Consequently, denoting with $\uds\in\uc^{N'}$ the dataset with $n'_i=N'$ and $n'_j=0$ for each $j\neq i$, a priori we have,
$$\lP(\realized'=\uds)=0,\quad \uP(\realized'=\uds)=1,$$
and in particular
$$\lP(\ux_1=x_i)=0,\quad \uP(\ux_1=x_i)=1.$$
It can be shown that other prior predictive probabilities
are not vacuous. For example, for $i\neq j$, we have
$$\uE(\theta_i\theta_j)=\frac{s}{4(s+1)}<\frac{1}{4}=\sup_{\theta\in\Theta} \theta_i\theta_j.$$

The IDM produces, for each possible observed data set $\uxs$,
non-vacuous posterior predictive probabilities for each possible
future data set (see \cite{Walley1996}). This means that our
previous theoretical limitations are avoided in some way. To explain
this result we consider two cases. We consider firstly an observed
data set $\uxs$ where we have observed at least two different
outcomes. Secondly, we consider a data set $\uxs$ formed exclusively
by outcomes of the same type, in other words, a data set of the type
$\uds$.

In the first case we have that $\P(\uxs\giv\theta)=\prod_{j=1}^k
\theta_j^{n_j}$ is equal to zero for $\theta=\uesi$ for each
$i\in\{1,\ldots ,k\}$. In fact, $\theta_i=1$ implies $\theta_j=0$
for each $j\neq i$ and there is at least one $j$ with $n_j>0$.
Therefore, the assumptions of Corollaries \ref{cor:fondamentale3}
and \ref{cor:fondamentale3bis} are not satisfied. And in fact the
IDM produces non-vacuous posterior predictive probabilities for each
data set that, a priori, has vacuous predictive probabilities. On
the other hand, all the datasets whose prior predictive probability
reaches its maximum in a relative frequency $\ufs\in\Theta$ such
that $\P(\uxs\giv\ufs)>0$, are characterized by non-vacuous prior
predictive probabilities.

The second case yields similar results. The only difference is that
$\P(\uds\giv\theta)=\theta_i^{N'}$ for a given $i\in\{1,\ldots
,k\}$. In this case $\P(\uxs\giv\uesi)=1>0$ and in fact, according
to Corollaries \ref{cor:fondamentale3} and
\ref{cor:fondamentale3bis}, we obtain
$$\uP(x_i\giv\uxs)=\uP(x_i)=1,$$
$$\uP(\ux'=\uds\giv\uxs)=\uP(\uds)=1,$$
and consequently, for each $j\neq i$ and each $\uys\neq\uds$,
$$\lP(x_j\giv\uxs)=\lP(x_j)=0,$$
$$\lP(\ux'=\uys\giv\uxs)=\lP(\uys)=0.$$
But, on the other hand, we obtain
$$\lP(x_i\giv\uxs)>0,\quad \lP(\ux'=\uds\giv\uxs)>0,$$
$$\uP(x_j\giv\uxs)<1,\quad \uP(\ux'=\uys\giv\uxs)<1,$$
and therefore the posterior predictive probabilities are not vacuous
for each possible future data set.\eoe
\end{exmp}

Yet, since the variables $(\ux_i)_{i\in\enne}$ are assumed to be
observable, the successful application of a near-ignorance prior in
the IDM is not helpful in addressing the doubts raised by our
theoretical results about the applicability of near-ignorance priors
in situations where the variables $(\ux_i)_{i\in\enne}$ are latent.

\section{Conclusions}\label{sec:conclusions}

In this paper we have proved a sufficient condition that prevents
learning about a latent categorical variable to take place under
prior near-ignorance about the data-generating process.

The condition holds as soon as the likelihood is strictly positive
(and continuous), and so is satisfied frequently, even in the
simplest settings. Taking into account that the considered framework
is very general and pervasive of statistical practice, we regard
this result as a form of substantial evidence against the possibility to
use prior near-ignorance in real statistical problems. Given that
complete prior ignorance is not compatible with learning, as it is well
known, we deduce that there is little hope to use any form of prior
ignorance to do objective-minded statistical inference in practice.

As a consequence, we suggest that future research efforts should be
directed to study and develop new forms of knowledge that are close
to near-ignorance but that do not coincide with it.
\section*{Acknowledgements}
This work was partially supported by Swiss NSF grants
200021-113820/1 (Alberto Piatti), 200020-109295/1 (Marco Zaffalon)
and 100012-105745/1 (Fabio Trojani).

\begin{appendix}

\section{Technical preliminaries}\label{techresults}

In this appendix we provide some technical results that are used to
prove the theorems in the paper. First of all, we introduce some
notation used in this appendix. Consider a sequence of probability
densities $\seqpn$ and a function $f$ defined on a set $\ste$. Then,
we use the notation,
$$\E_n(f):=\int_{\ste} f(\teta) p_n(\teta) d\teta,$$
$$\P_n(\widetilde{\ste}):=\int_{\widetilde{\ste}} p_n(\teta) d\teta,\quad \widetilde{\ste}\subseteq\ste.$$
In addition, for a given probability density $p$ on $\ste$,
$$\E_p(f):=\int_{\ste} f(\teta) p(\teta) d\teta,$$
$$\P_p(\widetilde{\ste}):=\int_{\widetilde{\ste}} p(\teta) d\teta,\quad \widetilde{\ste}\subseteq\ste.$$
Finally, with $\rightarrow$ we denote $\lim_{n\rightarrow\infty}$.

\begin{thm}\label{thm:marcus1}
Let $\ste\subset\erre^k$ be the closed $k$-dimensional simplex and
let $\seqpn$ be a sequence of probability densities defined on
$\ste$ w.r.t. the Lebesgue measure. Let $f\geq 0$ be a bounded
continuous function on $\ste$ and denote with $f_{\max}$ the
supremum of $f$ on $\ste$. For this function define the measurable
sets
\begin{equation}\ste_{\delta}=\{\teta\in\ste\giv f(\teta)\geq
f_{\max}-\delta\}.\label{deftetadelta}\end{equation} Assume that
$\seqpn$ concentrates on a maximum of $f$ for $n\rightarrow\infty$,
in the sense that
\begin{equation} \E_n(f)\rightarrow f_{\max},\label{convfmax}\end{equation} then, for all $\delta >0$, it
holds
$$\P_n(\ste_{\delta})\rightarrow 1.$$
\end{thm}

\begin{thm}\label{thm:marcus2}
Let $L(\teta)\geq 0$ be a bounded measurable function with
\begin{equation}
\lim_{\delta\rightarrow 0} \inf_{\teta\in\ste_{\delta}} L(\teta) =:
c>0, \label{liminfpos}
\end{equation}
under the same assumptions of Theorem \ref{thm:marcus1}. Then
\begin{equation}
\frac{\E_n(Lf)}{\E_n(L)}=\frac{\int_{\ste} f(\teta) L(\teta)
p_n(\teta) d\teta}{\int_{\ste} L(\teta) p_n(\teta)
d\teta}\rightarrow f_{\max}. \label{vacuousthm2}
\end{equation}
\end{thm}

\begin{rem}\label{rem:assumptionmarcus}
If $f$ has a unique maximum in $\teta=\teta_0$ and $L$ is a
function, continuous in an arbitrary small neighborhood of
$\theta=\theta_0$, such that $L(\teta_0)>0$, then (\ref{liminfpos})
is satisfied.
\end{rem}

\section{Corollaries to Theorem \ref{thm:fondamentale}}\label{app:corollaries}

The following Corollaries to Theorem \ref{thm:fondamentale} are
necessary to prove Corollary \ref{cor:fondamentale4}, and are useful
to understand more deeply the limiting results implied by the use of
near-ignorance priors with latent variables.

\begin{cor}\label{cor:fondamentale1}
Let $\uxs'$ and $\uss$ be given. Denote with
$\mathbf{f'}:=\left(\frac{n_1'}{N'},\ldots,\frac{n_k'}{N'}\right)\in\Theta$
the vector of relative frequencies of the dataset $\uxs'$. If
$\P(\uss\giv\theta)$ is continuous in an arbitrary small
neighborhood of $\theta=\mathbf{f'}$, $\P(\uss\giv\mathbf{f'})>0$
and $\mzero$  is such that
$$\uP(\mathbf{\uxs'})=\sup_{\theta\in\Theta} \left(\prod_{i=1}^k \theta_i^{n_i'}\right)
=\prod_{i=1}^k \left(\frac{n_i'}{N'}\right)^{n_i'},$$ then
$$\uP(\mathbf{\uxs'}\giv\uss)=\uP(\mathbf{\uxs'}).$$
\end{cor}

\begin{cor}\label{cor:fondamentale2}
Let $\uxs'$ and $\uss$ be given. If $\P(\uss\giv\theta)>0$ for each
$\theta\in\Theta$ with $\theta_i=0$ for at least one $i$ with
$n_i'>0$, and $\mzero$ is such that $\lP(\uxs')=0$, it
follows that
$$\lP(\uxs'\giv\uss)=\lP(\uxs')=0.$$
\end{cor}

\begin{cor}\label{cor:fondamentale3}
Let $\uss$ be given. Consider an arbitrary $x_i\in\uc$ and denote
with $\uesi$ the particular vector of chances with $\theta_i=1$ and
$\theta_j=0$ for each $j\neq i$. Suppose that $\mzero$ is such that,
a priori, $\uP(\ux_1=x_i):=\uE(\theta_i)=1$. Then,
if $\P(\uss\giv\uesi)>0$ and $\P(\uss\giv\theta)$ is continuous in a
neighborhood of $\theta=\uesi$, we have
\begin{equation}\label{eq:uppervacue}
\uP(\ux_{N+1}=x_i\giv\uss)=\uP(\ux_1=x_i)=1,
\end{equation}
and consequently,
\begin{equation}\label{eq:lowervacue}
\lP(\ux_{N+1}=x_j\giv\uss)=\lP(\ux_j=x_i)=0,
\end{equation}
for each $j\neq i$.
\end{cor}

\begin{cor}\label{cor:fondamentale3bis}
Let $\uss$ and $N'$ be given and consider an arbitrary $x_i\in\uc$.
Suppose that $\mzero$ is such that, a priori,
$\uP(\ux_1=x_i):=\uE(\theta_i)=1$. Denote with
$\uds\in\uc^{N'}$ the data set with $n_i=N'$ and $n_j=0$ for each
$j\neq i$. Then, if $\P(\uss\giv\uesi)>0$ and $\P(\uss\giv\theta)$ is
continuous in a neighborhood of $\theta=\uesi$, we have
$$\uP(\realized'=\uds\giv\uss)=1,$$
and consequently,
$$\lP(\realized'=\uys\giv\uss)=0,$$
for each $\uys\neq\uds$.
\end{cor}

\end{appendix}

\pagebreak[3]

\end{document}